\documentclass[12pt]{amsart}

\usepackage{amssymb,latexsym,amscd}

\usepackage{amsthm,amsmath}

\newcommand{\cO}{\mathcal{O}}

\newcommand{\cE}{\mathcal{E}}

\newcommand{\cM}{\mathcal{M}}

\newcommand{\cB}{\mathcal{B}}

\newcommand{\cL}{\mathcal{L}}

\newcommand{\Sym}{\mathrm{Sym}}
\newcommand{\Pic}{\mathrm{Pic}}

\newcommand{\lra}{\longrightarrow}

\newcommand{\ra}{\rightarrow}

\newcommand{\la}{\lambda}

\newcommand{\PP}{\mathbb{P}}
\newcommand{\Sbb}{\mathbb{S}}

\newcommand{\CC}{\mathbb{C}}

\newcommand{\End}{\mathrm{End}}

\newcommand{\Ext}{\mathrm{Ext}}

\newcommand{\im}{\mathrm{im}}
\newcommand{\rk}{\mathrm{rk}}

\def\map#1{\ \smash{\mathop{\longrightarrow}\limits^{#1}}\ }

\theoremstyle{plain}

\newtheorem{thm}{Theorem}[section]

\newtheorem{lem}[thm]{Lemma}

\newtheorem{prop}[thm]{Proposition}

\newtheorem{cor}[thm]{Corollary}

\begin{document}

\title[]{Rank four vector bundles without theta divisor over a curve of genus two}

\author{Christian Pauly}

\address{D\'epartement de Math\'ematiques \\ Universit\'e de Montpellier II - Case Courrier 051 \\ Place Eug\`ene Bataillon \\ 34095 Montpellier Cedex 5 \\ France}
\email{pauly@math.univ-montp2.fr}




\begin{abstract}
We show that the locus of stable rank four vector bundles without theta divisor over a smooth projective curve of genus two is in canonical bijection with the set of theta-characteristics. We give several descriptions of these bundles and compute the degree of the rational theta map. 
\end{abstract}

\maketitle


\section{Introduction}
Let $C$ be a complex smooth projective curve of genus $2$ and let $\cM_r$  denote the
coarse moduli space parametrizing semi-stable rank-$r$ vector bundles with trivial 
determinant over the curve $C$. Let $C \cong \Theta \subset \Pic^1(C)$ be the Riemann theta
divisor in the degree $1$ component of the Picard variety of $C$. For any $E \in \cM_r$
we consider the locus
$$ \theta(E) = \{ L \in \Pic^1(C)  \ | \ h^0(C, L \otimes E ) > 0 \}, $$
which is either a curve linearly equivalent to $r\Theta$ or $\theta(E) = \Pic^1(C)$, in
which case we say that $E$ has no theta divisor. We obtain thus a rational map, the so-called
theta map
$$\theta : \cM_r \dasharrow |r \Theta|,$$
between varieties having the same dimension $r^2 -1$. We denote by $\cB_r$ the closed subvariety
of $\cM_r$ parametrizing semi-stable bundles without theta divisor. It is known \cite{R} that $\cB_2= \cB_3 = \emptyset$
and that $\cB_r \not= \emptyset$ for $r \geq 4$.

\bigskip

Is was recently shown that $\theta$ is generically finite; see \cite{B2} Theorem A. Moreover the cases of low ranks $r$ have been studied in the past:
if $r=2$ the theta map is an isomorphism $\cM_2 \cong \PP^3$ \cite{NR}  and if $r=3$ the theta
map realizes $\cM_3$ as a double covering of $\PP^8$ ramified along a 
sextic hypersurface \cite{O}.

\bigskip

In this note we study the next case $r=4$ and give a complete description of the locus $\cB_4$.
Our main result is the following 

\begin{thm}
Let $C$ be a curve of genus $2$.
\begin{enumerate}
\item The locus $\cB_4$ is of dimension $0$, reduced and of cardinality $16$.
\item There exists a canonical bijection between $\cB_4$ and the set of 
theta-characteristics of $C$. Let $E_\kappa \in \cB_4$ denote the stable vector
bundle associated with the theta-characteristic $\kappa$. Then
$$ \Lambda^2 E_\kappa = \bigoplus_{\alpha \in S(\kappa)} \alpha,
\qquad \Sym^2 E_\kappa = \bigoplus_{\alpha \in J[2] \setminus S(\kappa)} \alpha,
$$
where $S(\kappa)$ is the set of $2$-torsion line bundles $\alpha \in J[2]$
such that $\kappa \alpha \in \Theta \subset \Pic^1(C)$.

\item If $\kappa$ is odd, then $E_\kappa$ is a symplectic bundle. If $\kappa$ is even,
then $E_\kappa$ is an orthogonal bundle with non-trivial Stiefel-Whitney class.

\item The $16$ vector bundles $E_\kappa$ are invariant under the tensor product 
with the group $J[2]$. 
\end{enumerate}
\end{thm}

The $16$ vector bundles $E_\kappa$ already appeared in Raynaud's paper \cite{R} as
Fourier-Mukai transforms and were further studied in \cite{Hi} and \cite{He}
--- see section 2.2. We note that Theorem 1.1 completes the main result of \cite{Hi}
which describes the restriction of $\cB_4$ to {\em symplectic} rank-$4$ bundles. 
The method of this paper is different and is partially based on \cite{P}.

\bigskip

As an application of Theorem 1.1 we obtain the degree of the theta map for $r=4$. We 
refer to \cite{BV} for a geometric interpretation of the general fiber of $\theta$ in terms of
certain irreducible components of a Brill-Noether locus of the curve $\theta(E) \subset
\Pic^1(C)$.  

\begin{cor}
The degree of the rational theta map $\theta :  \cM_4 \dasharrow |4\Theta|$ equals $30$.
\end{cor}

\bigskip

{\em Acknowledgements:} I am grateful to George Hitching and Olivier Serman for useful discussions.

\bigskip

{\em Notations:} If $E$ is a vector bundle over $C$, we will write $H^i(E)$ for $H^i(C,E)$
and $h^i(E)$ for $\dim H^i(C,E)$. We denote the slope of $E$ by $\mu(E) := 
\frac{\deg E}{\rk E}$, the canonical bundle over $C$ by $K$ and the degree $d$
component of the Picard variety of $C$ by $\Pic^d(C)$. We denote by $J:= \Pic^0(C)$ the
Jacobian of $C$ and by $J[n]$ its group of $n$-torsion points. The divisor 
$\Theta_\kappa \subset J$ is the translate of the Riemann theta divisor $C \cong \Theta \subset
\Pic^1(C)$ by a theta-characteristic $\kappa$. 
The line bundle $\cO_J(2\Theta_\kappa)$ does not depend on $\kappa$ and will be denoted
by $\cO_J(2\Theta)$.

\bigskip

\section{Proof of Theorem 1.1}

\subsection{The $16$ vector bundles $E_\kappa$}

We first show that the set-theoretical support of $\cB_4$ consists of 
$16$ stable vector bundles $E_\kappa$, which are canonically labelled by the
theta-characteristics of $C$.

\bigskip

We note that $\cB_4 \not= \emptyset$ by \cite{R}, see also \cite{P} Theorem 1.1. 
We consider a vector
bundle $\cE \in \cB_4$. Since $\cB_2 = \cB_3 = \emptyset$, we deduce that $\cE$ is stable.
We introduce $\cE' = \cE^* \otimes K$. Then $\mu(\cE') = 2$ and since 
$\cE \in \cB_4$, we obtain that $h^0(\cE' \otimes \la^{-1}) = h^1(\cE \otimes
\la) = h^0(\cE \otimes \la) > 0$ for any $\la \in \Pic^1(C)$. In particular
for any $x \in C$ we have $h^0(\cE' \otimes \cO_C(-x)) > 0$.  On the other hand
stability of $\cE$ implies that $h^0(\cE) = h^1(\cE') =  0$. Hence $h^0(\cE') = 4$
by Riemann-Roch. Thus we obtain that the evaluation map of global sections
$$ \cO_C \otimes H^0(\cE') \map{ev} \cE' $$
is not of maximal rank. Let us denote by $I := \im \ ev$ the subsheaf of $\cE'$
given by the image of $ev$. Then clearly $h^0(I) = 4$. The cases $\rk \ I \leq 2$ are 
easily ruled out using stability of $\cE'$. Hence we conclude that $\rk \ I = 3$. 
We then consider the natural exact sequence 
\begin{equation} \label{es3}
0 \lra L^{-1} \lra \cO_C \otimes H^0(\cE') \map{ev} I \lra 0,
\end{equation}
where $L$ is the line bundle such that $L^{-1} := \ker \ ev$.

\begin{prop}
We have $h^0(I^*)= 0$.
\end{prop}

\begin{proof}
Suppose on the contrary that there exists a non-zero map $I \ra \cO_C$. Its kernel
$S \subset I$ is a rank-$2$ subsheaf of $\cE'$ and by stability of $\cE'$  we obtain 
$\mu(S) < \mu(\cE') = 2$, hence $\deg S \leq 3$. Moreover $h^0(S) \geq h^0(I) -1 = 3$.

\bigskip

Assume that $\deg S = 3$. Then $S$ is stable and $S$ can be written as an extension
$$ 0 \lra \mu \lra S \lra \nu \lra 0,$$
with $\deg \mu = 1$ and $\deg \nu = 2$. The condition $h^0(S) \geq 3$ then
implies that $\mu = \cO_C(x)$ for some $x \in C$, $\nu = K$ and that the 
extension has to be split, i.e., $S = K \oplus \cO_C(x)$. This contradicts stability
of $S$.

\bigskip

The assumption $\deg S \leq 2$ similarly leads to a contradiction. We leave the 
details to the reader. \end{proof}

Now we take the cohomology of the dual of the exact sequence \eqref{es3} and we 
obtain --- using $h^0(I^*) = 0$ --- an inclusion $H^0(\cE')^* \subset H^0(L)$. 
Hence $h^0(L) \geq 4$, which implies $\deg L \geq 5$. On the other hand $\deg L = 
\deg I$ and by stability of $\cE'$, we have $\mu(I) < 2$, i.e., $\deg L \leq 5$.
So we can conclude that $\deg L = 5$, that $H^0(\cE')^* = H^0(L)$ and that
$I = E_L$, where $E_L$ is the {\em evaluation bundle} associated to $L$ defined
by the exact sequence
\begin{equation} \label{esEL}
 0 \lra E^*_L \lra H^0(L) \otimes \cO_C \map{ev} L \lra 0. 
\end{equation}

\bigskip

Moreover the subsheaf $E_L \subset \cE'$ is of maximal degree, hence $E_L$ is a subbundle
of $\cE'$ and we have an exact sequence 
\begin{equation} \label{ese}
 0 \lra E_L \lra \cE' \lra K^4 L^{-1} \lra 0,
\end{equation}
with extension class $e \in \Ext^1(K^4 L^{-1}, E_L) = H^1(E_L \otimes K^{-4}L) = 
H^0(E^*_L \otimes K^5 L^{-1})^*$. Using Riemann-Roch and stability of $E_L$ (see e.g. \cite{Bu})
one shows that
$$h^0(E^*_L \otimes K^5 L^{-1})  = 7, \qquad h^0(E^*_L \otimes K^5 L^{-1}(-x)) = 4,  \qquad 
h^0(E^*_L \otimes K^5 L^{-1}(-x-y)) = 1$$ for {\em general} points $x,y \in C$. In that case we 
denote by $\mu_{x,y} \in \PP H^0(E^*_L \otimes K^5 L^{-1})$ the point determined by the
$1$-dimensional subspace $H^0(E^*_L \otimes K^5 L^{-1}(-x-y))$. We also denote by 
$$\Sbb \subset \PP H^0(E^*_L \otimes K^5 L^{-1})$$ 
the linear span of the points $\mu_{x,y}$ when
$x$ and $y$ vary in $C$ and by $H_e \subset \PP H^0(E^*_L \otimes K^5 L^{-1})$ the
hyperplane determined by the non-zero class $e$.

\bigskip

Tensoring the sequence \eqref{ese}  with $K^{-4}L(x+y)$  and taking cohomology one shows that $\mu_{x,y} \in H_e$
if and only if $h^0(\cE' \otimes K^{-4}L(x+y)) > 0$. Since we assume $\cE \in \cB_4$,
we obtain 
$$\Sbb \subset H_e.$$
We consider a {\em general} point $x \in C$ such that $h^0(E^*_L \otimes K^5 L^{-1}(-x)) = 4$
and denote for simplicity 
$$ A := E^*_L \otimes K^5 L^{-1}(-x). $$
Then $A$ is stable with $\mu(A) = \frac{7}{3}$. We consider the evaluation map of global
sections 
$$ ev_A : \cO_C \otimes H^0(A) \lra A$$
and consider the set $S_A$ of points $p \in C$ for which  $(ev_A)_p$ is not surjective, i.e.
$$ S_A = \{ p \in C \ | \ h^0(A(-p)) \geq 2 \}.$$
Then we have the following 

\begin{lem}
We assume that $x$ is general.
\begin{enumerate}
\item If $L^2 \not= K^5$, then the set $S_A$ consists of the $2$ distinct points $p_1,p_2$ determined by
the relation $\cO_C(p_1 + p_2) = K^4 L^{-1}(-x)$.
\item If $L^2 = K^5$, then the set $S_A$ consists of the $2$ distinct points $p_1,p_2$ introduced in (1) and
the conjugate $\sigma(x)$ of $x$ under the hyperelliptic involution $\sigma$. 
\end{enumerate}
\end{lem} 

\begin{proof}
Given a point $p \in C$, we tensorize the exact sequence \eqref{esEL} with 
$K^5 L^{-1} (-x-p)$ and take cohomology:
$$ 0 \lra H^0(A(-p)) \lra H^0(L) \otimes H^0(K^5 L^{-1}(-x-p)) \lra 
H^0(K^5(-x-p)) \lra \cdots $$
We note that $h^0(K^5 L^{-1}(-x-p)) = 2$. We distinguish two cases.

\bigskip
\noindent
(a) The pencil $|K^5 L^{-1} (-x-p)|$ has a base-point, i.e. there exists a point $q \in C$
such that  $K^5 L^{-1} (-x-p) = K(q)$, or equivalently $K^4L^{-1}(-x) = \cO_C(p+q)$. Since
$x$ is general, we have  $h^0(K^4L^{-1}(-x)) = 1$, which determines $p$ and $q$, i.e., 
$\{ p,q \} = \{ p_1, p_2 \}$. In this case $|K^5 L^{-1} (-x-p)| = |K(q)| = |K|$ and
$h^0(A(-p)) = h^0(K^{-1}L) = 2$. This shows that $p_1,p_2 \in S_A$.

\bigskip
\noindent
(b) The pencil $|K^5 L^{-1} (-x-p)|$ is base-point-free. By the base-point-free-pencil-trick,
we have $H^0(A(-p)) \cong H^0(L^2 K^{-5} (x+p))$. Since $\deg L^2 K^{-5} (x+p) = 2$, we have
$h^0(L^2 K^{-5} (x+p)) = 2$ if and only if $L^2 K^{-5} (x+p) = K$, or equivalently
$\cO_C(p) = K^6 L^{-2}(-x)$. If $K^6 L^{-2} \not= K$, then for general $x \in C$ the 
line bundle
$K^6 L^{-2}(-x)$ is not of the form $\cO_C(p)$. If $K^6 L^{-2} = K$, then for any $x \in C$,
$K^6L^{-2}(-x) = \cO_C(\sigma(x))$, which implies that $\sigma(x) \in S_A$.

\bigskip

This shows the lemma.
\end{proof}

\begin{prop}
If $L^2 \not= K^5$, then $\Sbb =  \PP H^0(E^*_L \otimes K^5 L^{-1})$.
\end{prop}

\begin{proof}
We consider a general point $x \in C$ and the rank-$3$ bundle $A$. Let $B \subset A$ 
denote the subsheaf given by the image of $ev_A$. By Lemma 2.2 (1) we have 
$\deg B = \deg A - 2 = 5$. Moreover $H^0(B) = H^0(A)$ and there is an
exact sequence 
\begin{equation} \label{es4}
0 \lra M^{-1} \lra \cO_C \otimes H^0(B) \map{ev_A} B \lra 0,
\end{equation}
with $M \in \Pic^5(C)$. It follows that the rational map 
$$ \phi_x : C \dasharrow \PP H^0(B) = \PP H^0(A) = \PP^3, \qquad y \mapsto \mu_{x,y}$$
factorizes through
$$ C \map{\varphi_M} |M|^* \lra \PP H^0(B),$$
where $\varphi_M$ is the morphism given by the linear system $|M|$
and the second map is linear and identifies with the projectivization of the dual
of $\delta$, which is given by the long exact sequence obtained from \eqref{es4} by
dualizing and taking cohomology:
$$ 0 \lra H^0(B^*) \lra H^0(B)^* \map{\delta} H^0(M) \lra H^1(B^*) \lra \cdots$$
We obtain that the linear span of $\im \  \phi_x$ is non-degenerate if and only if
$h^0(B^*) = 0$.

\bigskip

We now show that $h^0(B^*)= 0$. Suppose on the contrary that there exists a non-zero map
$B \ra \cO_C$. Its kernel $S \subset B$ is a rank-$2$ subsheaf of $A$ with $\deg S 
\geq \deg B = 5$, hence $\mu(S) \geq \frac{5}{2}$, which contradicts stability of $A$ ---
recall that $\mu(A) = \frac{7}{3}$.

\bigskip

This shows that $\im \ \phi_x$ spans $\PP H^0(A) \subset \PP H^0(E_L^* \otimes 
K^5 L^{-1})$ for general $x \in C$. We now take $2$ general points $x,x' \in C$ and 
deduce from $\dim H^0(A) \cap H^0(A') = \dim H^0(E^*_L \otimes K^5 L^{-1}(-x-x')) = 1$
that the linear span of the union $\PP H^0(A) \cap \PP H^0(A')$ equals the full space
$\PP H^0(E_L^* \otimes K^5 L^{-1})$. This shows the proposition.
\end{proof}   

We deduce from the proposition that the line bundle $L$ satisfies the relation
$L^2 = K^5$, i.e. 
$$ L = K^2 \kappa$$
for some theta-characteristic $\kappa$ of $C$. In that case we note that
$H^0(E_L^* \otimes K^5 L^{-1})$ equals $H^0(E_L^* \otimes L)$ and we can consider the 
exact sequence 
$$ 0 \lra H^0(E_L^* \otimes L) \lra H^0(L) \otimes H^0(L) \map{\mu} H^0(L^2) \lra 0,$$
obtained from \eqref{esEL} by tensoring with $L$ and taking cohomology. We also note that there is a
natural inclusion  $\Lambda^2 H^0(L) \subset H^0(E^*_L \otimes L)$, see e.g. \cite{P} section 2.1. 
More precisely we can show 

\begin{prop}
The linear span $\Sbb$ equals
$$\Sbb = \PP \Lambda^2 H^0(L) \subset \PP H^0(E^*_L \otimes L).$$
\end{prop}

\begin{proof}
Using the standard exact sequences and the base-point-free-pencil-trick, one
easily works out that for general points $x,y \in C$
$$ \mu_{x,y} = \PP \Lambda^2 H^0(L(-x-y)) \subset \PP \Lambda^2 H^0(L)
\subset \PP H^0(E^*_L \otimes L).$$
This implies that $\Sbb \subset \PP \Lambda^2 H^0(L)$. In order to show 
equality one chooses $4$ general points $x_i \in C$ such that their images 
$C \ra |L|^* = \PP^3$ linearly span the $\PP^3$. We denote by $s_i \in H^0(L)$ 
the global section vanishing on the points $x_j$ for $j \not= i$ and not
vanishing on $x_i$. Then one checks that for any choice of the indices $i,j,k,l$ 
such that $\{i,j,k,l\} = \{1,2,3,4 \}$ one has $s_i \wedge s_j = \mu_{x_k,x_l}$.
Since the $6$ tensors $s_i \wedge s_j$ are a basis of $\Lambda^2 H^0(L)$, we obtain
equality.
\end{proof} 

The hyperplane $\Sbb = \PP \Lambda^2 H^0(L) \subset \PP H^0(E^*_L \otimes L)$ determines
a unique (up to a scalar) non-zero extension class $e \in H^0(E^*_L \otimes L)^*$ by 
$\Sbb = H_e$, which in turn determines a unique stable vector bundle $\cE \in \cB_4$, which 
we will denote by $E_\kappa$.

\bigskip

This shows that $\cB_4$ is of dimension $0$ and of cardinality $16$.

\subsection{The Raynaud bundles}

In this subsection we recall the construction of the Raynaud bundles introduced in 
\cite{R} as Fourier-Mukai transforms. We refer to \cite{Hi} section 9.2 for the
details and the proofs.

\bigskip

The rank-$4$ vector bundle $\cO_J(2\Theta) \otimes H^0(J, \cO_J(2 \Theta))^*$ over $J$ admits a
canonical $J[2]$-lineari- zation and descends therefore under the duplication map
$[2]: J \ra J$, i.e., there exists a rank-$4$ vector bundle $M$ over $J$ such that
$$[2]^* M \cong \cO_J(2\Theta) \otimes H^0(J, \cO_J(2 \Theta))^*.$$

\begin{prop}
For any theta-characteristic $\kappa$ of $C$ there exists an isomorphism
$$\xi_\kappa : M \map{\sim} M^* \otimes \cO_J(\Theta_\kappa).$$
Moreover if $\kappa$ is even (resp. odd), then $\xi_\kappa$ is symmetric
(resp. skew-symmetric).
\end{prop}

Let $\gamma_\kappa : C \ra J$ be the Abel-Jacobi map defined by $\gamma_\kappa(p) = 
\kappa^{-1}(p)$. We define the Raynaud bundle 
$$ R_\kappa := \gamma_\kappa^* M \otimes \kappa^{-1}.$$
Then by \cite{R} the bundle $R_\kappa \in \cB_4$. Since $\gamma_\kappa^* \cO_J(\Theta_\kappa) 
= K$ we see that the isomorphism $\xi_\kappa$ induces an orthogonal (resp. symplectic)
structure on the bundle $R_\kappa$, if $\kappa$ is even (resp. odd). In particular the
bundle $R_\kappa$ is self-dual, i.e., $R_\kappa = R_\kappa^*$. The pull-back
$\gamma^*_\kappa(\xi_\kappa')$ for a theta-characteristic $\kappa' = \kappa \alpha$
with $\alpha \in J[2]$ gives an isomorphism
$$ R_\kappa \map{\sim} R^*_\kappa \otimes \alpha,$$
hence a non-zero section in $H^0(\Lambda^2 R_\kappa \otimes \alpha)$ (resp. 
$H^0(\Sym^2 R_\kappa \otimes \alpha)$) if $h^0(\kappa \alpha) = 1$ (resp. 
$h^0(\kappa \alpha) = 0$). We deduce that there are isomorphisms
\begin{equation} \label{isomRkappa}
\Lambda^2 R_\kappa = \bigoplus_{\alpha \in S(\kappa)} \alpha,
\qquad \Sym^2 R_\kappa = \bigoplus_{\alpha \in J[2] \setminus S(\kappa)} \alpha.
\end{equation}
In particular the $16$ bundles $R_\kappa$ are non-isomorphic. Each $R_\kappa$ is 
invariant under tensor product with $J[2]$. The isomorphisms \eqref{isomRkappa}
can be used to prove the relation 
\begin{equation} \label{invaj4}
R_\kappa \otimes \beta = R_{\kappa \beta^2}, \qquad \forall \beta \in J[4]. 
\end{equation}  

\subsection{Symplectic and orthogonal bundles}

In this subsection we give a third construction of the bundles in $\cB_4$ as
symplectic and orthogonal extension bundles. Let $\kappa$ be a theta-characteristic.

\bigskip

If $\kappa$ is odd, then $\kappa = \cO_C(w)$ for some Weierstrass point $w \in C$. 
The construction outlined in \cite{P} section 2.2 gives a unique symplectic bundle $\cE_e \in \cB_4$ 
with $e \in H^1(\Sym^2 G)_+$. We denote this bundle by $V_\kappa$.

\bigskip

If $\kappa$ is even, there is an analogue construction, which we briefly outline for 
the convenience of the reader. The proofs are similar to those given in \cite{Hi}.
Using the Atiyah-Bott-fixed-point formula one observes that among all non-trivial
extensions
$$ 0 \lra \kappa^{-1} \lra G \lra \cO_C \lra 0,$$
there are $2$ extensions (up to scalar), 
which are $\sigma$-invariant. We take one of them. Then
any non-zero class $e \in H^1(\Lambda^2 G) = H^1(\kappa^{-1})$ determines
an orthogonal bundle $\cE_e$, which fits in the exact sequence
\begin{equation} \label{esG}
 0 \lra G \lra \cE_e \lra G^* \lra 0.
\end{equation}
The composite map  
$$ D_G : \PP H^1 (\Lambda^2 G) \lra \cM_4 \map{\theta} |4\Theta|, \qquad e \mapsto
\theta(\cE_e)$$
is the projectivization of a linear map
$$ \widetilde{D_G} : H^1(\Lambda^2 G) \lra H^0(\Pic^1(C), 4 \Theta). $$
Moreover $\im \ \widetilde{D_G} \subset H^0(\Pic^1(C), 4 \Theta)_-$, which can be
seen as follows. By \cite{Se} Thm 2 the second Stiefel-Whitney class $w_2(\cE_e)$ of an
orthogonal bundle $\cE_e$ is given by the parity of $h^0(\cE_e \otimes \kappa')$ for
any theta-characteristic $\kappa'$. This parity can be computed by taking the cohomology
of the exact sequence \eqref{esG} tensorized with $\kappa'$ and taking into account that
the coboundary map is skew-symmetric. One obtains that $w_2(\cE_e) \not= 0$ and one can
conclude the above-mentioned inclusion by \cite{B3} Lemma 1.4.

\bigskip

We now observe that by the Atiyah-Bott-fixed-point-formula $h^1(\Lambda^2 G)_+ =
h^1(\Lambda^2 G)_- = 1$. By the argument given in \cite{P} section 2.2 we conclude that one
of the two eigenspaces $H^1(\Lambda^2 G)_\pm$ is contained in the kernel $\ker \ 
\widetilde{D_G}$. We denote the corresponding bundle $\cE_e$ by $V_\kappa \in \cB_4$.

\subsection{Three descriptions of the same bundle}

\begin{prop}
For any theta-characteristic $\kappa$ the three bundles $E_\kappa$, $R_\kappa$ and $V_\kappa$ coincide.
\end{prop}

\begin{proof}
If $\kappa$ is odd, this was worked out in detail in \cite{Hi} section 8 and Theorem 29. If 
$\kappa$ is even, the proofs are similar.
\end{proof} 

This proposition shows all assertions of Theorem 1.1 except reducedness of $\cB_4$.

\bigskip
 
I am grateful to Olivier Serman for giving me the following fourth description of the  bundle
$E_\kappa$ for an  even  theta-characteristic $\kappa$. We recall that an even theta-characteristic
$\kappa$ corresponds to a partition of the set of six Weierstrass points of $C$ into two subsets
of three points, which we denote by $\{w_1,w_2,w_3 \}$ and $\{w_4,w_5,w_6 \}$. With this notation we have 

\begin{prop}
Let $\kappa$ be an even theta-characteristic. We denote by $A_\kappa$  (resp. $B_\kappa$) the unique stable
rank-$2$ bundle with determinant $\kappa$ and which contains the four $2$-torsion line bundles
$\cO_C$, $\cO_C(w_1 - w_2)$, $\cO_C(w_1 - w_3)$ and $\cO_C(w_2 - w_3)$ (resp.                  
$\cO_C$, $\cO_C(w_4 - w_5)$, $\cO_C(w_4 - w_6)$ and $\cO_C(w_5 - w_6)$). Then the orthogonal rank-4 vector bundle $E_\kappa$ is 
isomorphic to 
$$ \mathrm{Hom}(A_\kappa, B_\kappa) $$
equipped with the quadratic form given by the determinant. 
\end{prop}

We refer to \cite{S} section 5.5 for the proof.

\bigskip

\subsection{Reducedness of $\cB_4$}

We denote by $\cL$ the determinant line bundle over the moduli space $\cM_4$ and
recall that the set $\cB_4$ can be identified with the base locus of the linear
system $|\cL|$. This endows the set $\cB_4$ with a natural scheme-structure. 

\bigskip

We start with a description of the space of global sections $H^0(\cM_4, \cL)$.

\begin{prop}
For any theta-characteristic $\kappa$ there is a section $s_\kappa \in H^0(\cM_4, \cL)$ with
zero divisor
$$ \Delta_\kappa :=\mathrm{Zero}(s_\kappa) 
= \{ E \in \cM_4 \ | \ h^0(\Lambda^2 E \otimes \kappa) > 0 \}.$$
The $16$ sections $s_\kappa$ form a basis of $H^0(\cM_4, \cL)$.
\end{prop}

\begin{proof}
The Dynkin index of the second fundamental representation $\rho : \mathfrak{sl}_4(\CC) \ra 
\End(\Lambda^2 \CC^4)$ equals $2$ (see e.g. \cite{LS} Proposition 2.6). Moreover the bundle
$\Lambda^2 E \otimes \kappa$ admits a $K$-valued non-degenerate quadratic form, which allows
to construct the Pfaffian divisor $s_\kappa$, which is a section of $\cL$ (see \cite{LS}).
The space $H^0(\cM_4, \cL)$ is a representation of level $2$ of the Heisenberg group $Heis(2)$,
which is a central extension of $J[2]$ by $\CC^*$. One can work out that the sections
$s_\kappa$ generate the $16$ one-dimensional character spaces for the $Heis(2)$-action on
$H^0(\cM_4, \cL)$. This shows that the sections $s_\kappa$ are linearly independent.
\end{proof}

Since $E_\kappa \in \cB_4$, we have $E_\kappa \in \Delta_{\kappa'}$ for any theta-characteristic
$\kappa'$. By the deformation theory of determinant and Pfaffian divisors (see e.g. \cite{L},
\cite{LS}) the  point $E_\kappa \in \cM_4$ is a smooth point of the divisor $\Delta_{\kappa'}
\subset \cM_4$ if and only if the following two conditions hold
\begin{enumerate}
\item $h^0(\Lambda^2 E_\kappa \otimes \kappa') = 2$,
\item the natural linear form
$$ \Phi_{\kappa'} : T_{E_\kappa} \cM_4 = H^1(\End_0(E_\kappa)) \lra \Lambda^2 
H^0(\Lambda^2 E_\kappa \otimes \kappa')^*$$
is non-zero.
\end{enumerate}
Moreover if these two conditions holds, then $T_{E_\kappa} \Delta_{\kappa'} = \ker 
\Phi_{\kappa'}$.
The map $\Phi_{\kappa'}$ is built up as follows: the exceptional isomorphism
of Lie algebras $\mathfrak{sl}_4 \cong \mathfrak{so}_6$ induces a natural 
vector bundle isomorphism 
\begin{equation} \label{isomEkappa}
\End_0(E_\kappa) \map{\sim} \Lambda^2 ( \Lambda^2 E_\kappa ).
\end{equation}
Then $\Phi_{\kappa'}$ is the dual of the linear map given by the wedge product of global sections
$$\Lambda^2 H^0(\Lambda^2 E_\kappa \otimes \kappa') \lra H^0 (\Lambda^2 ( \Lambda^2 E_\kappa )
\otimes K) = H^0 (\End_0(E_\kappa) \otimes K).$$

\begin{prop}
The $0$-dimensional scheme $\cB_4$ is reduced.
\end{prop}

\begin{proof}
Since $E_\kappa$ is a smooth point of $\cM_4$ and $\dim T_{E_\kappa} \cM_4 = 15$, it is sufficient
to show that for any theta-characteristic $\kappa' \not= \kappa$ the divisor 
$\Delta_{\kappa'}$ is smooth at $E_\kappa$ and that the $15$ hyperplanes 
$\ker \Phi_{\kappa'} \subset T_{E_\kappa} \cM_4$ are linearly independent: using the 
isomorphism \eqref{isomRkappa} we obtain that for $\kappa' \not= \kappa$
$$ h^0(\Lambda^2 E_\kappa \otimes \kappa') = \sharp S(\kappa) \cap S(\kappa') = 2 $$
and using the isomorphism \eqref{isomEkappa} we obtain that
$$ \End_0(E_\kappa) = \bigoplus_{\alpha \in J[2] \setminus \{ 0 \} }  \alpha.$$
On the other hand one easily sees that if $\gamma, \delta \in J[2]$ are the two $2$-torsion
points in the intersection $S(\kappa) \cap S(\kappa')$, then $\kappa' = \kappa
\gamma \delta$, hence $\Lambda^2 H^0(\Lambda^2 E_\kappa \otimes \kappa') \cong
H^0(K \gamma \delta)$. This implies that the linear form
$$ \Phi_{\kappa'} : \bigoplus_{\alpha \in J[2] \setminus \{ 0 \} }  H^1(\alpha)
\lra H^0(K \gamma \delta)^* = H^1(\beta) $$
is projection onto the direct summand $H^1(\beta)$, where $\beta = \kappa^{-1} \kappa'
\in J[2]$. This description of the linear forms $\Phi_{\kappa'}$ clearly shows  that
they are non-zero and linearly independent.
\end{proof}

This completes the proof of Theorem 1.1. 

\bigskip

\section{Proof of Corollary 1.2}

Since by Theorem 1.1 $\cB_4$ is a reduced $0$-dimensional scheme of length $16$, the 
degree of the theta map $\theta$ is given by the formula
$$ \deg \theta + 16  = c_{15},$$
where $\frac{c_{15}}{15!}$ is the leading coefficient of the Hilbert polynomial
$$ P(n) = \chi(\cM_4, \cL^{n}) = \frac{c_{15}}{15!} n^{15} + \text{lower degree terms}.$$
In order to compute the polynomial $P$ we write 
\begin{equation} \label{exprP}
P(X) = \sum_{k=0}^{15} \alpha_k Q_k(X), \qquad \text{with} \qquad 
Q_k(X) = \frac{1}{k!} (X+7)(X+6) \cdots (X+8-k) 
\end{equation}
and $Q_0(X) = 1$. Note that $\deg Q_k = k$ and that $c_{15} = \alpha_{15}$. The 
canonical bundle of $\cM_4$ equals $\cL^{-8}$. By the Grauert-Riemenschneider
vanishing theorem we obtain that $h^i(\cM_4, \cL^{n}) = 0$ for any $i \geq 1$ and
$n \geq -7$. Hence $P(n) = h^0(\cM_4, \cL^n)$  for $n \geq -7$. Moreover
$P(n) = 0$ for $n = -7, -6, \ldots , -1$ and $P(0) = 1$. The values $P(n)$ for
$n=1,2, \ldots, 8$ can be computed by the Verlinde formula and with the use of MAPLE.
They are given in the following table.

\bigskip

\begin{center}
\begin{tabular}{|c||c|c|c|c|c|c|c|c|}
\hline 
 $n$ & 1 & 2 & 3 & 4 & 5 & 6 & 7 & 8 \\
\hline
 $P(n)$ & 16 & 140 & 896 & 4680 & 21024 & 83628 & 300080 & 984539 \\
\hline 
\end{tabular}
\end{center}

\bigskip

Using the expression \eqref{exprP} of $P$ one straightforwardly deduces the coefficients
$\alpha_k$ by increasing induction on $k$: $\alpha_k = 0$ for $k= 0,1, \ldots, 6$
and the values $\alpha_k$ for $k= 7, \ldots , 15$ are given in the following table.

\bigskip

\begin{center}
\begin{tabular}{|c||c|c|c|c|c|c|c|c|c|}
\hline 
 $k$ & 7 & 8 & 9 & 10 & 11 & 12 & 13 & 14 & 15 \\
\hline
 $\alpha_k$ & 1 & 8 & 32 & 96 & 214 & 328 & 324 & 184 & 46 \\
\hline 
\end{tabular}
\end{center}

\bigskip

Hence $\deg \theta = \alpha_{15} -16 = 30$.

\end{document}